# An efficient technique for fractional modified Boussinesq and approximate long wave equations

P. Veeresha[1], D. G. Prakasha[1,*], M. A. Qurashi[2], D. Baleanu[3,4]

Abstract

In this paper, an efficient technique is employed to study the modified Boussinesq and approximate long wave equations of the Caputo fractional time derivative, namely $q$-homotopy analysis transform method. These equations are playing a vital rule in describing the properties of shallow water waves through distinct dispersion relation. The convergence analysis and error analysis has been presented in the present investigation for the future scheme. We illustrate two examples to demonstrate the leverage and effectiveness of the proposed scheme, and the error analysis has been discussed to verify the accuracy. The numerical simulation has been conducted to ensure the exactness of the future technique. The obtained numerical and graphical results are divulge, the proposed scheme is computationally very accurate and straightforward to study and find the solution for fractional coupled nonlinear complex phenomena arised in science and technology.

**Keywords:** Laplace transform; fractional approximate long wave equations; fractional modified Boussinesq equations; $q$-homotopy analysis transform method.

## 1. Introduction

Fractional calculus (FC) was firstly put forward by L'Hopital. FC has garnered a lot of devotion and appreciation during the last few decades, due to their capability to provide an exact description for various non-linear complex phenomena. The differential systems with fractional order have lately gained popularity in developing procedure of models and investigation of dynamical systems. The fractional calculus is the generalization of the traditional calculus having nonlocal and genetic consequence in the material properties.

The fundamental properties of fractional calculus are described by many researchers [1-4]. FC plays a vital role and act as an essential tool in analysing and solving problems situated in diverse areas of science and technology like, fluid and continuum mechanics [5], chaos theory [6], biotechnology [7], electrodynamics [8], and many other fields [9-11]. The solution for differential equations having order arbitrary describing above phenomena plays a pivotal part in labelling the behaviour of complex problems arises in nature.

*Corresponding Author:
[1]Department of Mathematics, Faculty of Science & Technology, Karnatak University, Dharwad - 580003, India
[2]Department of Mathematics, College of Science, King Saud University, Riyadh - 11495, Saudi Arabia
[3]Department of Mathematics, Cankaya University, Balgat, Ankara - 06530, Turkey
[4]Institute of Space Sciences, Magurele, Bucharest - 077125, Romania
E-mail addresses: viru0913@gmail.com (P. Veeresha), prakashadg@gmail.com, dgprakasha@kud.ac.in (D. G. Prakasha), maysaa@ksu.edu.sa (M. A. Qurashi) dumitru@cankaya.edu.tr (D. Baleanu).

In twentieth century, Whitham [12], Broer [13] and Kaup [14] studied the equations, which elucidate the propagation of shallow water waves having distinct dispersion relation, called as Whitham-Broer-Kaup (WBK) equations. Consider the coupled WBK equations of fractional order [15]:

$$\begin{cases} D_t^\alpha u + uu_x + v_x + bu_{xx} = 0, \\ D_t^\alpha v + uv_x + vu_x + au_{xxx} - bv_{xx} = 0, \end{cases} \quad 0 < \alpha \leq 1, \tag{1}$$

where $u = u(x,t)$ is the horizontal velocity and $v = v(x,t)$ be the height that deviating from equilibrium position of the liquid. Here, $\alpha$ is the order of the time-fractional derivative. Further, $a$ and $b$ are constants which are representing distinct diffusion powers i.e., if $a = 1$ and $b = 0$, then Eq. (1) becomes modified Boussinesq equation. Similarly, for $a = 0$ and $b = 1$, the system signifies conventional long wave equation. These equations arise in hydrodynamics to illustrate the propagation of waves in dissipative and nonlinear media, and they are advisable for problems arise in the leakage of water in porous subsurface stratum and widely used in ocean and coastal engineering. Moreover, Eq. (1) is the foundation of numerous models utilize to portray the unconfined subsurface like, drainage and groundwater flow problems.

Last thirty years has been the testimony for discovery of plenty of new schemes to solve non-linear fractional differential equations in parallel to the developments of new computational algorithms with symbolic programming. In connection with this, Liao proposed a technique called homotopy analysis method (HAM) and which is based on construction of a homotopy which continuously deforms an initial guess approximation to the exact solution of the given problem [16]. It does not require any discretization, linearization and perturbation. But, it requires more computation and computer memory to solve nonlinear problems arise in complex phenomena. Hence, it necessitates a mixture of transformation algorithm to overcome with these confines.

In the present investigation, we employ $q$-HATM to find an approximated analytical solution for coupled modified Boussinesq and approximate long wave equations of fractional order. These equations are studied by several authors via distinct techniques like, ADM [17], VIM [18], CFRDTM [15], LADM [19] and other techniques [20-22]. These cited techniques require huge computation and has highly complicated procedure to solve coupled systems. The proposed method is a modified technique and which is elegant blend of $q$-HAM with Laplace transform. Hence, it does not require discretization, linearization or perturbation and in additionally it will decrease huge mathematical computations, more computer memory and is free from obtaining difficult integrations, polynomials, physical parameters. The future technique has many sturdy properties including straight forward solution procedure, promising large convergence region, and moreover free from any assumption, discretization and perturbation. It is worth revealing that the Laplace transform with semi-analytical techniques requires less C.P.U time to evaluate solution for nonlinear complex models and phenomena arised in science and technology. The $q$-HATM solution involves with two auxiliary parameters $\hbar$ and $n$, which helps us to adjust and control the convergence of the solution, which quickly tends to the analytical solution in a small acceptable region. It is worth to mention that, the proposed scheme can decrease the computation of the time and work as compared with other traditional techniques while maintaining the great efficiency.

Recently, due to consistency and efficacy of $q$-HATM it has been eminently aided by many researchers to analyse various kinds of nonlinear problem like, authors in [23] analysed and find the approximated analytical solution for fractional model of vibration equation and elucidate the exactness of $q$-HATM, in [24] the solution for fractional Drinfeld–Sokolov–Wilson equation has been investigated by the aid of proposed scheme, the cancer chemotherapy effect model with fractional order has been analysed by the authors in [25], H. Bulut and his co-authors analyse HIV infection of CD4+T lymphocyte cells of fractional model [26], the efficiency of considered algorithm is presented while finding the solution for fractional telegraph equation equation in [27], authors in [28] analysed the model of Lienard's equation, and many others has study and find the solution for many complex problems arised in related fields of science [29-34]. Motivated by these investigations, we find the approximated analytical solution for coupled DEs describing the shallow water waves through distinct dispersion relation.

## 2. Preliminaries

In this segment, we present basic definitions and notions which will be used in the present frame work:

**Definition 1** Let $f(t) \in C_\mu$ ($\mu \geq -1$) be a function. Then Riemann-Liouville fractional integral of $f(t)$ with order $\alpha > 0$ presented as [4]

$$J^\alpha f(t) = \frac{1}{\Gamma(\alpha)} \int_0^t (t-x)^{\alpha-1} f(\vartheta) dx,$$
$$J^0 f(t) = f(t). \quad (2)$$

**Definition 2** The Caputo fractional derivative of $f \in C_{-1}^n$ is defined as

$$D_t^\alpha f(t) = \begin{cases} \frac{d^n f(t)}{dt^n}, & \alpha = n \in \mathbb{N}, \\ \frac{1}{\Gamma(n-\alpha)} \int_0^t (t-x)^{n-\alpha-1} f^{(n)}(x) dx, & n-1 < \alpha < n, n \in \mathbb{N}. \end{cases} \quad (3)$$

**Definition 3** The Laplace transform (LT) of a Caputo fractional derivative $D_t^\alpha f(t)$ is represented as [2, 3]

$$L[D_t^\alpha f(t)] = s^\alpha F(s) - \sum_{r=0}^{n-1} s^{\alpha-r-1} f^{(r)}(0^+), (n-1 < \alpha \leq n), \quad (4)$$

where $F(s)$ is represent the LT of $f(t)$.

## 3. Basic idea of $q$-HATM

Consider a nonlinear differential equation of arbitrary order

$$D_t^\alpha \mathcal{U}(x,t) + R\,\mathcal{U}(x,t) + N\,\mathcal{U}(x,t) = f(x,t), \quad n-1 < \alpha \leq n, \quad (5)$$

where $R$ is the bounded linear differential operator in $x$ and $t$, (i.e., for a number $\varepsilon > 0$ we have $\|R\mathcal{U}\| \leq \varepsilon\|\mathcal{U}\|$), $N$ specifies the nonlinear differential operator and Lipschitz continuous with $\mu > 0$ satisfying $|N\mathcal{U} - N\mathcal{V}| \leq \mu|\mathcal{U} - \mathcal{V}|$, and $f(x,t)$ denotes the source term. On applying $LT$ to Eq. (5), we have

$$s^\alpha L[\mathcal{U}(x,t)] - \sum_{k=0}^{n-1} s^{\alpha-k-1}\mathcal{U}^{(k)}(x,0) + L[R\mathcal{U}(x,t)] + L[N\mathcal{U}(x,t)] = L[f(x,t)]. \quad (6)$$

On simplification, the Eq. (6) reduces to

$$L[\mathcal{U}(x,t)] - \frac{1}{s^\alpha}\sum_{k=0}^{n-1} s^{\alpha-k-1}\mathcal{U}^k(x,0) + \frac{1}{s^\alpha}\{L[R\mathcal{U}(x,t)] + L[N\mathcal{U}(x,t)] - L[f(x,t)]\} = 0. \quad (7)$$

According to homotopy analysis method [16], the nonlinear operator for is real function $\varphi(x,t;q)$ defined as

$$N[\varphi(x,t;q)] = L[\varphi(x,t;q)] - \frac{1}{s^\alpha}\sum_{k=0}^{n-1} s^{\alpha-k-1}\varphi^{(k)}(x,t;q)(0^+)$$
$$+ \frac{1}{s^\alpha}\{L[R\varphi(x,t;q)] + L[N\varphi(x,t;q)] - L[f(x,t)]\}, \quad \forall q \in \left[0, \frac{1}{n}\right] \quad (8)$$

The homotopy constructed for $H(x,t)$ as shown below:

$$(1 - nq)L[\varphi(x,t;q) - \mathcal{U}_0(x,t)] = \hbar q H(x,t) N[\varphi(x,t;q)], \quad (9)$$

where $L$ is symbolized LT, $\hbar \neq 0$ is an auxiliary parameter, $\mathcal{U}_0(x,t)$ is initial guess and $\varphi(x,t;q)$ is an unknown function. For $q = 0$ and $q = \frac{1}{n}$, the following results are respectively holds true

$$\varphi(x,t;0) = \mathcal{U}_0(x,t), \quad \varphi\left(x,t;\frac{1}{n}\right) = \mathcal{U}(x,t), \quad (10)$$

Thus, by increasing $q$ from 0 to $\frac{1}{n}$, the solution $\varphi(x,t;q)$ converge from $\mathcal{U}_0(x,t)$ to $\mathcal{U}(x,t)$. Now by applying Taylor theorem, the function $\varphi(x,t;q)$ is expanding in series form near to $q$ [35], we have

$$\varphi(x,t;q) = \mathcal{U}_0(x,t) + \sum_{m=1}^{\infty} \mathcal{U}_m(x,t) q^m, \quad (11)$$

where

$$\mathcal{U}_m(x,t) = \frac{1}{m!}\frac{\partial^m \varphi(x,t;q)}{\partial q^m}\Big|_{q=0}. \quad (12)$$

On choosing the initial guess $\mathcal{U}_0(x,t)$, the auxiliary parameter $n$, the auxiliary linear operator, $\hbar$ and $H(x,t)$, the series (11) converges at $q = \frac{1}{n}$. Later, it provides solutions for Eq. (5) and which is of the form

$$\mathcal{U}(x,t) = \mathcal{U}_0(x,t) + \sum_{m=1}^{\infty} \mathcal{U}_m(x,t)\left(\frac{1}{n}\right)^m. \tag{13}$$

Now, differentiating Eq. (9) $m$-times in terms of $q$ and then multiplied by $\frac{1}{m!}$ and then taking $q = 0$, it yields

$$L[\mathcal{U}_m(x,t) - k_m \mathcal{U}_{m-1}(x,t)] = \hbar H(x,t)\Re_m(\vec{\mathcal{U}}_{m-1}), \tag{14}$$

where

$$\vec{\mathcal{U}}_m = \{\mathcal{U}_0(x,t), \mathcal{U}_1(x,t), \ldots, \mathcal{U}_m(x,t)\}. \tag{15}$$

On employing inverse LT for Eq. (14), one can get

$$\mathcal{U}_m(x,t) = k_m \mathcal{U}_{m-1}(x,t) + \hbar L^{-1}[H(x,t)\Re_m(\vec{\mathcal{U}}_{m-1})], \tag{16}$$

where

$$\Re_m(\vec{\mathcal{U}}_{m-1}) = L[\mathcal{U}_{m-1}(x,t)] - \left(1 - \frac{k_m}{n}\right)\left(\sum_{k=0}^{n-1} s^{\alpha-k-1}\mathcal{U}^{(k)}(x,0) + \frac{1}{s^\alpha}L[f(x,t)]\right) + \frac{1}{s^\alpha}L[R\mathcal{U}_{m-1} + \mathcal{H}_{m-1}], \tag{17}$$

and

$$k_m = \begin{cases} 0, & m \leq 1, \\ n, & m > 1. \end{cases} \tag{18}$$

In Eq. (17), $\mathcal{H}_m$ denotes homotopy polynomial and defined as

$$\mathcal{H}_m = \frac{1}{m!}\left[\frac{\partial^m \varphi(x,y,t;q)}{\partial q^m}\right]_{q=0} \text{ and } \varphi(x,y,t;q) = \varphi_0 + q\varphi_1 + q^2\varphi_2 + \cdots. \tag{19}$$

By Eqs. (16) and (17), we have

$$\mathcal{U}_m(x,,t) = (k_m + \hbar)\mathcal{U}_{m-1}(x,t) - \left(1 - \frac{k_m}{n}\right)L^{-1}(\sum_{k=0}^{n-1} s^{\alpha-k-1}\mathcal{U}^{(k)}(x,0) + \frac{1}{s^\alpha}L[f(x,t)]) + \hbar L^{-1}\left[\frac{1}{s^\alpha}L[R\mathcal{U}_{m-1} + \mathcal{H}_{m-1}]\right]. \tag{20}$$

Lastly, on simplifying Eq. (20) we get the iterative terms of $\mathcal{U}_m(x,t)$. The series solution of $q$-HATM is presented by

$$\mathcal{U}(x,t) = \sum_{m=0}^{\infty} \mathcal{U}_m(x,t). \tag{21}$$

## 4. Convergence analysis of $q$-HATM solution

**Theorem 1.** (Uniqueness theorem)

The solution for the non-linear fractional differential equation (5) obtained by $q$-HATM is unique for every $\beta \in (0,1)$, where $\beta = (n + \hbar) + \hbar(\varepsilon + \mu)T$.

**Proof:** For Eq. (5), the solution is defined by

$$\mathcal{U}(x,t) = \sum_{m=0}^{\infty} \mathcal{U}_m(x,t),$$

where
$$\mathcal{U}_m(x,t) = (k_m + \hbar)\mathcal{U}_{m-1}(x,t) - \left(1 - \frac{k_m}{n}\right)L^{-1}\left(\sum_{k=0}^{n-1} s^{\alpha-k-1}\mathcal{U}^{(k)}(x,0)\right.$$
$$\left. + \frac{1}{s^\alpha}L[f(x,t)]\right) + \hbar L^{-1}\left[\frac{1}{s^\alpha}L[R\mathcal{U}_{m-1} + \mathcal{H}_{m-1}]\right]. \tag{22}$$

Suppose $\mathcal{U}$ and $\mathcal{U}^\blacksquare$ are the two solutions of Eq. (4), then it is sufficient to show $\mathcal{U} = \mathcal{U}^\blacksquare$ to prove the theorem, Now by Eq. (21), we obtain

$$|\mathcal{U} - \mathcal{U}^\blacksquare| = \left|(n + \hbar)(\mathcal{U} - \mathcal{U}^\blacksquare) + \hbar L^{-1}\left(\frac{1}{s^\alpha}L\big(N(\mathcal{U} - \mathcal{U}^\blacksquare) + R(\mathcal{U} - \mathcal{U}^\blacksquare)\big)\right)\right|,$$

then by using the convolution theorem for $LT$, we have

$$|\mathcal{U} - \mathcal{U}^\blacksquare| \leq (n + \hbar)|\mathcal{U} - \mathcal{U}^\blacksquare| + \hbar \int_0^t (|N(\mathcal{U} - \mathcal{U}^\blacksquare)| + |R(\mathcal{U} - \mathcal{U}^\blacksquare)|) \frac{(t-\xi)^\alpha}{\Gamma(\alpha+1)} d\xi$$

$$\leq (n + \hbar)|\mathcal{U} - \mathcal{U}^\blacksquare| + \hbar \int_0^t (\varepsilon|(\mathcal{U} - \mathcal{U}^\blacksquare)| + \mu|(\mathcal{U} - \mathcal{U}^\blacksquare)|) \frac{(t-\xi)^\alpha}{\Gamma(\alpha+1)} d\xi.$$

By the aid of integral mean value theorem, the above equation reduces to

$$|\mathcal{U} - \mathcal{U}^\blacksquare| \leq (n + \hbar)|\mathcal{U} - \mathcal{U}^\blacksquare| + \hbar(\varepsilon|(\mathcal{U} - \mathcal{U}^\blacksquare)| + \mu|(\mathcal{U} - \mathcal{U}^\blacksquare)|)T.$$

Here $\beta = (n + \hbar) + \hbar(\varepsilon + \mu)T$, thus

$$|\mathcal{U} - \mathcal{U}^\blacksquare| \leq \beta|\mathcal{U} - \mathcal{U}^\blacksquare| \Rightarrow (1 - \beta)|\mathcal{U} - \mathcal{U}^\blacksquare| \leq 0$$

Since $0 < \beta < 1$, then $\mathcal{U} - \mathcal{U}^\blacksquare = 0 \Rightarrow \mathcal{U} = \mathcal{U}^\blacksquare$.

Hence, the solution for Eq. (5) is unique. $\square$

**Theorem 2**. (Convergence theorem)

Let $X$ be a Banach space and $F: X \to X$ be a non-linear mapping. Assume that

$$\|F(\mathcal{U}) - F(\mathcal{V})\| \leq \beta \|\mathcal{U} - \mathcal{V}\|, \ \forall \ a, b \in X,$$

then $F$ has a fixed point in view of Banach fixed point theory [36]. Moreover, for the arbitrary selection of $a_0, b_0 \in X$, the sequence generated by the $q$-HATM converges to fixed point of $F$ and

$$\|\mathcal{U}_m - \mathcal{U}_n\| \leq \frac{\beta^n}{1-\beta} \|\mathcal{U}_1 - \mathcal{U}_0\|, \ \forall \ a, b \in X.$$

**Proof:** For all continuous functions, let us consider Banach space $(C[I], \|\cdot\|)$ on $I$ with norm is given by $\|g(\lambda)\| = \max_{\lambda \in I}|g(\lambda)|$. First, we prove that $\{\mathcal{U}_n\}$ is Cauchy sequence in $X$.

Now consider,

$$\|\mathcal{U}_m - \mathcal{U}_n\| = \max_{\lambda \in I}|\mathcal{U}_m - \mathcal{U}_n|$$

$$= \max_{\lambda \in I}\left|(n + \hbar)(\mathcal{U}_{m-1} - \mathcal{U}_{n-1}) + \hbar L^{-1}\left(\frac{1}{s^\alpha}L(N(\mathcal{U}_{m-1} - \mathcal{U}_{n-1}) + R(\mathcal{U}_{m-1} - \mathcal{U}_{n-1}))\right)\right|$$

$$\leq \max_{\lambda \in I}[(n + \hbar)|(\mathcal{U}_{m-1} - \mathcal{U}_{n-1})| + \hbar L^{-1}(\frac{1}{s^\alpha}L(N|\mathcal{U}_{m-1} - \mathcal{U}_{n-1}|$$

$$+ R(|\mathcal{U}_{m-1} - \mathcal{U}_{n-1}|))). \tag{23}$$

By the convolution theorem for $LT$, the Eq. (22) becomes

$$\|\mathcal{U}_m - \mathcal{U}_n\| \leq \max_{\lambda \in I}[(n+\hbar)|(\mathcal{U}_{m-1} - \mathcal{U}_{n-1})| + \hbar \int_0^t (|N(\mathcal{U}_{m-1} - \mathcal{U}_{n-1})| + |R(\mathcal{U}_{m-1} - \mathcal{U}_{n-1})|) \frac{(t-\xi)^\alpha}{\Gamma(\alpha+1)} d\xi$$

$$\leq \max_{\lambda \in I}[(n+\hbar)|(\mathcal{U}_{m-1} - \mathcal{U}_{n-1})| + \hbar \int_0^t (|N(\mathcal{U}_{m-1} - \mathcal{U}_{n-1})| + |R(\mathcal{U}_{m-1} - \mathcal{U}_{n-1})|) \frac{(t-\xi)^\alpha}{\Gamma(\alpha+1)} d\xi. \tag{24}$$

By the aid of integral mean value theorem, the above inequality reduces to

$$\|\mathcal{U}_m - \mathcal{U}_n\| \leq \max_{\lambda \in I}[(n+\hbar)|(\mathcal{U}_{m-1} - \mathcal{U}_{n-1})| + \hbar(\varepsilon|\mathcal{U}_{m-1} - \mathcal{U}_{n-1}| + \mu|\mathcal{U}_{m-1} - \mathcal{U}_{n-1}|)T]$$

$$\leq \beta \|\mathcal{U}_{m-1} - \mathcal{U}_{n-1}\|.$$

For $m = n+1$, one can get

$$\|\mathcal{U}_{n+1} - \mathcal{U}_n\| \leq \beta \|\mathcal{U}_n - \mathcal{U}_{n-1}\| \leq \beta^2 \|\mathcal{U}_{n-1} - \mathcal{U}_{n-2}\| \leq \beta^3 \|\mathcal{U}_{n-2} - \mathcal{U}_{n-3}\| \leq \cdots$$

$$\leq \beta^n \|\mathcal{U}_1 - \mathcal{U}_0\|.$$

In view of triangular inequality, we have

$$\|\mathcal{U}_m - \mathcal{U}_n\| \leq \|\mathcal{U}_{n+1} - \mathcal{U}_n\| + \|\mathcal{U}_{n+2} - \mathcal{U}_{n+1}\| + \cdots + \|\mathcal{U}_m - \mathcal{U}_{m-1}\|$$

$$\leq [\beta^n + \beta^{n+1} + \cdots + \beta^{m-1}]\|\mathcal{U}_1 - \mathcal{U}_0\| = \beta^n[1 + \beta + \cdots + \beta^{m-n-1}]\|\mathcal{U}_1 - \mathcal{U}_0\|$$

$$\leq \beta^n \left[\frac{1-\beta^{m-n-1}}{1-\beta}\right] \|\mathcal{U}_1 - \mathcal{U}_0\|.$$

Clearly, $1 - \beta^{m-n-1} < 1$ (since $0 < \beta < 1$). Therefore, the above inequality becomes

$$\|\mathcal{U}_m - \mathcal{U}_n\| \leq \frac{\beta^n}{1-\beta} \|\mathcal{U}_1 - \mathcal{U}_0\|. \tag{25}$$

But $\|\mathcal{U}_1 - \mathcal{U}_0\| < \infty$, consequently as $m \to \infty$ then $\|\mathcal{U}_m - \mathcal{U}_n\| \to 0$.

It provides $\{\mathcal{U}_n\}$ is Cauchy sequence in $C[I]$, and every Cauchy sequence is convergent sequence. Hence $\{\mathcal{U}_n\}$ is convergent sequence. □

## 5. Error analysis of proposed algorithm

The error analysis of proposed scheme obtained with help of $q$-HATM is presented in this segment

**Theorem 3:** If we can obtain a real number $0 < \rho < 1$ fulfilling $\|\mathcal{U}_{m+1}(x,t)\| \leq \rho \|\mathcal{U}_m(x,t)\|$ for all $m$. Moreover, if the truncated series $\sum_{m=0}^{i} \mathcal{U}_m(x,t)$ is used as approximate solution of $\mathcal{U}(x,t)$, then maximum absolute truncated error can be obtained by

$$\left\|\mathcal{U}(x,t) - \sum_{m=0}^{i}\mathcal{U}_m(x,t)\right\| \leq \frac{\rho^{i+1}}{1-\rho}\|\mathcal{U}_0(x,t)\|$$

**Proof:** We have

$$\left\|\mathcal{U}(x,t) - \sum_{m=0}^{i}\mathcal{U}_m(x,t)\right\| = \left\|\sum_{m=i+1}^{\infty}\mathcal{U}_m(x,t)\right\|$$

$$\leq \sum_{m=l+1}^{\infty}\|\mathcal{U}_m(x,t)\|$$

$$\leq \sum_{m=i+1}^{\infty}\rho^m\|\mathcal{U}_0(x,t)\|$$

$$\leq \rho^{i+1}[1 + \rho^1 + \rho^2 + \cdots]\|\mathcal{U}_0(x,t)\|$$

$$\leq \frac{\rho^{i+1}}{1-\rho}\|\mathcal{U}_0(x,t)\|, \qquad \square$$

which proves our required result.

## 6. Illustrative examples

Here, we consider two coupled examples to present the efficiency and applicability $q$-HATM.

**Example 6.1.**

Consider the modified Boussinesq (MB) equations of fractional order [15, 17, 20]:

$$\begin{cases} D_t^\alpha u = -u\dfrac{\partial u}{\partial x} - \dfrac{\partial v}{\partial x}, \\ D_t^\alpha v = -u\dfrac{\partial v}{\partial x} - v\dfrac{\partial u}{\partial x} - \dfrac{\partial^3 u}{\partial x^3}, \end{cases} \quad 0 < \alpha \leq 1, \tag{26}$$

with initial conditions

$$u(x,0) = \omega - 2\ell \coth[\ell(x+c)], \qquad v(x,0) = -2\ell^2 \operatorname{csch}^2[\ell(x+c)]. \tag{27}$$

The exact solution for classical oreder MB equations is

$$u(x,t) = \omega - 2\ell \coth[\ell(x+c-\omega t)], \quad v(x,t) = -2\ell^2 \operatorname{csch}^2[\ell(x+c-\omega t)].$$

By performing LT on Eq. (26) and using Eq. (27), we have

$$\begin{aligned} L[u(x,t)] - \tfrac{1}{s}(\omega - 2\ell \coth[\ell(x+c)]) + \tfrac{1}{s^\alpha} L\left\{ u\tfrac{\partial u}{\partial x} + \tfrac{\partial v}{\partial x} \right\} &= 0, \\ L[v(x,t)] - \tfrac{1}{s}(-2\ell^2 \operatorname{csch}^2[\ell(x+c)]) + \tfrac{1}{s^\beta} L\left\{ u\tfrac{\partial v}{\partial x} + v\tfrac{\partial u}{\partial x} + \tfrac{\partial^3 u}{\partial x^3} \right\} &= 0. \end{aligned} \tag{28}$$

Define the non-linear operators as

$$\begin{aligned} N^1[\varphi_1(x,t;q), \varphi_2(x,t;q)] &= L[\varphi_1(x,t;q)] - \tfrac{1}{s}(\omega - 2\ell \coth[\ell(x+c)]) \\ &\quad + \tfrac{1}{s^\alpha} L\left\{ \varphi_1(x,t;q) \tfrac{\partial \varphi_1(x,t;q)}{\partial x} + \tfrac{\partial \varphi_2(x,t;q)}{\partial y} \right\}, \\ N^2[\varphi_1(x,t;q), \varphi_2(x,t;q)] &= L[\varphi_2(x,t;q)] - \tfrac{1}{s}(-2\ell^2 \operatorname{csch}^2[\ell(x+c)]) \\ &\quad + \tfrac{1}{s^\beta} L\left\{ \varphi_1(x,t;q) \tfrac{\partial \varphi_2(x,t;q)}{\partial x} + \varphi_2(x,t;q) \tfrac{\partial \varphi_1(x,t;q)}{\partial y} \right. \\ &\quad \left. + \tfrac{\partial^3 \varphi_1(x,t;q)}{\partial x^3} \right\}. \end{aligned} \tag{29}$$

By applying the proposed numerical scheme, $m$-th order deformation equation for $H(x,t)=1$, is given as

$$\begin{aligned} L[u_m(x,t) - k_m u_{m-1}(x,t)] &= h\mathfrak{R}_{1,m}[\vec{u}_{m-1}, \vec{v}_{m-1}], \\ L[v_m(x,t) - k_m v_{m-1}(x,t)] &= h\mathfrak{R}_{2,m}[\vec{u}_{m-1}, \vec{v}_{m-1}], \end{aligned} \tag{30}$$

where

$$\mathfrak{R}_{1,m}[\vec{u}_{m-1}, \vec{v}_{m-1}] = L[u_{m-1}(x,t)] - \left(1 - \tfrac{k_m}{n}\right) \tfrac{1}{s}(\omega - 2\ell \coth[\ell(x+c)]) \tag{31}$$

$$+\frac{1}{s^\alpha} L\left\{\sum_{i=0}^{m-1} u_i \frac{\partial u_{m-1-i}}{\partial x} + \frac{\partial v_{m-1}}{\partial x}\right\},$$

$$\Re_{2,m}[\vec{u}_{m-1},\vec{v}_{m-1}] = L[v_{m-1}(x,t)] - \left(1 - \frac{k_m}{n}\right)\frac{1}{s}(-2\ell^2\ csch^2[\ell(x+c)])$$

$$+\frac{1}{s^\beta} L\left\{\sum_{i=0}^{m-1} u_i \frac{\partial v_{m-1-i}}{\partial x} + \sum_{i=0}^{m-1} v_i \frac{\partial u_{m-1-i}}{\partial x} + \frac{\partial^3 u_{m-1}}{\partial x^3}\right\}.$$

By employing inverse LT on Eq. (30), we get

$$u_m(x,t) = k_m u_{m-1}(x,t) + \hbar L^{-1}\{\Re_{1,m}[\vec{u}_{m-1},\vec{v}_{m-1}]\},$$
$$v_m(x,t) = k_m v_{m-1}(x,t) + \hbar L^{-1}\{\Re_{2,m}[\vec{u}_{m-1},\vec{v}_{m-1}]\}. \tag{32}$$

On solving above equation, we have

$$u_0(x,t) = \omega - 2\ell\coth[\ell(x+c)],\ v_0(x,t) = -2\ell^2\ csch^2[\ell(x+c)],$$

$$u_1(x,t) = \frac{2\hbar\ell^2\omega\ csch^2[\ell(x+c)]t^\alpha}{\Gamma[\alpha+1]},\ v_1(x,t) = \frac{4\hbar\ell^3\omega\ coth[\ell(x+c)]\ csch^2[\ell(x+c)]t^\alpha}{\Gamma[\alpha+1]},$$

$$u_2(x,t) = \frac{2(n+\hbar)\hbar\ell^2\omega\ csch^2[\ell(x+c)]t^\alpha}{\Gamma[\alpha+1]} - \frac{4\hbar^2\ell^3\omega^2\ coth[\ell(x+c)]csch^2[\ell(x+c)]t^{2\alpha}}{\Gamma[2\alpha+1]},$$

$$v_2(x,t) = \frac{4(n+\hbar)\hbar\ell^3\omega\ coth[\ell(x+c)]\ csch^2[\ell(x+c)]t^\alpha}{\Gamma[\alpha+1]} - \frac{4\hbar^2\ell^4\omega^2(2+cosh[2\ell(x+c)])\ csch^4[\ell(x+c)]t^{2\alpha}}{\Gamma[2\alpha+1]},$$

$$u_3(x,t) = \frac{2(\hbar+n)^2\hbar\ell^2\omega\ csch^2[\ell(x+c)]\ t^\alpha}{\Gamma[1+\alpha]} - \frac{8(\hbar+n)\hbar^2\ell^3\omega^2\ coth[\ell(x+c)]\ csch^2[\ell(x+c)]\ t^{2\alpha}}{\Gamma[1+2\alpha]}$$

$$+\frac{2\hbar^3\ell^4\omega^2\ csch^5[\ell(x+c)]t^{3\alpha}}{\Gamma[\alpha+1]^2\ \Gamma[3\alpha+1]}(-4l\ cosh[\ell(x+c)]\ \Gamma[1+2\alpha]$$

$$+\Gamma[1+\alpha]^2\big(8\ell\ cosh[\ell(x+c)] + \omega((3\ sinh[\ell(x+c)] + sinh[3\ell(x+c)])\big)\big),$$

$$v_3(x,t) = \frac{4(\hbar+n)^2\hbar\ell^3\omega\ coth[\ell(x+c)]\ csch^2[\ell(x+c)]\ t^\alpha}{\Gamma[\alpha+1]} - \frac{8(\hbar+n)\hbar^2\ell^5\omega(2+cosh[2\ell(x+c)])\ csch^4[\ell(x+c)]t^{2\alpha}}{\Gamma[2\alpha+1]}$$

$$+\frac{2\hbar^3\ell^5\omega^2\ csch^6[\ell(x+c)]t^{3\alpha}}{\Gamma[\alpha+1]^2\ \Gamma[3\alpha+1]}(-4l(3 + 2\ cosh[2\ell(x+c)])\Gamma[2\alpha+1]$$

$$+\Gamma[\alpha+1]^2(24\ell + 16\ell\ cosh[2\ell(x+c)] + 10\omega\ sinh[2\ell(x+c)] + \omega\ sinh[4\ell(x+c)],$$

In this way, the rest term can be obtained. Then for system of Eq. (26), the $q$-HATM series solution of presented as follows

$$u(x,t) = u_0(x,t) + \sum_{m=1}^{\infty} u_m(x,t)\left(\frac{1}{n}\right)^m,$$
$$v(x,t) = v_0(x,t) + \sum_{m=1}^{\infty} v_m(x,t)\left(\frac{1}{n}\right)^m. \tag{33}$$

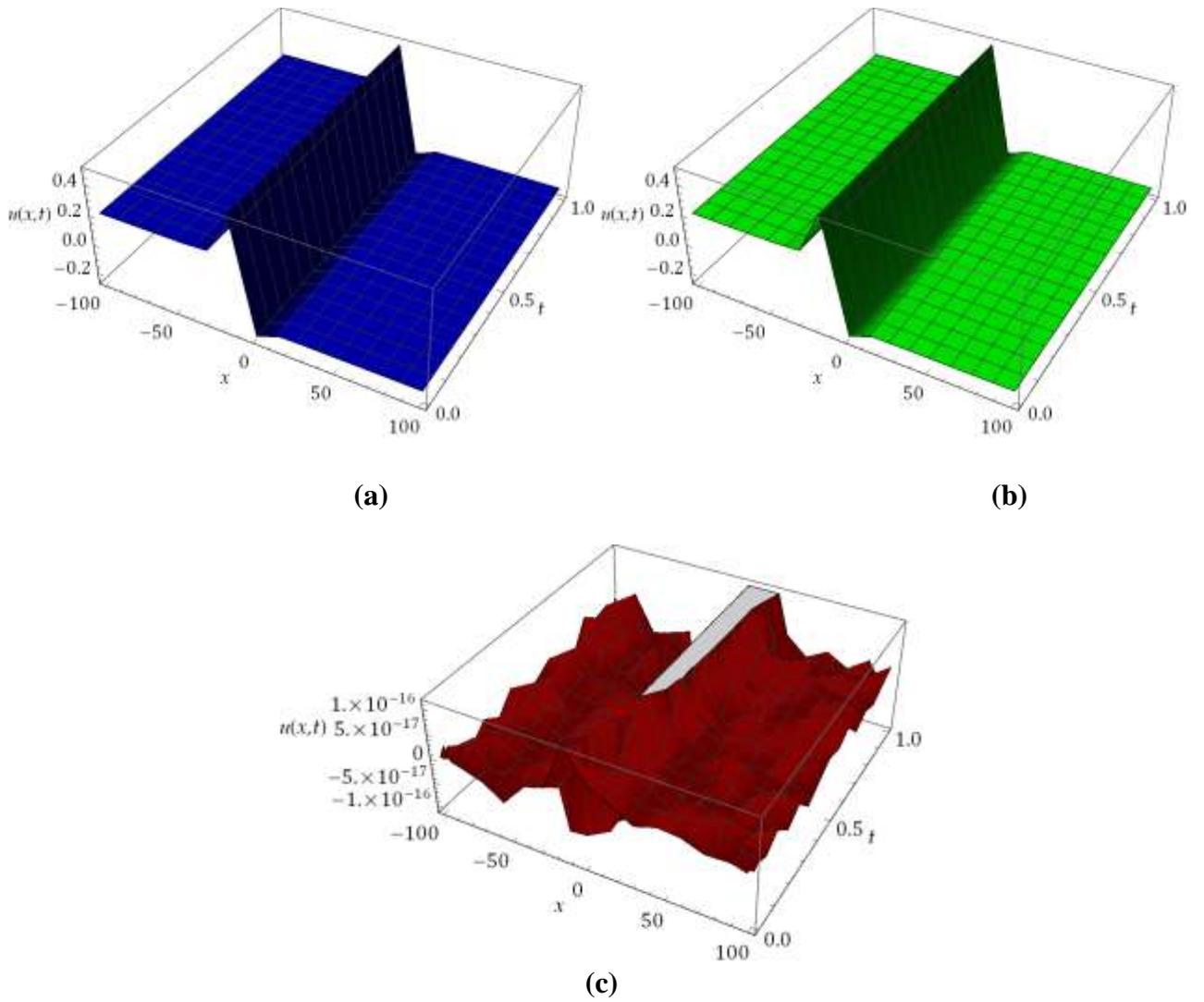

**Fig. 1**. Behaviour of **(a)** Approximate solution **(b)** Exact solution **(c)** Absolute error = $|u_{exa.} - u_{app.}|$ for Ex. 6.1 at $\omega = 0.005, \ell = 0.1, c = 10, n = 1, \alpha = 1$ and $\hbar = -1$.

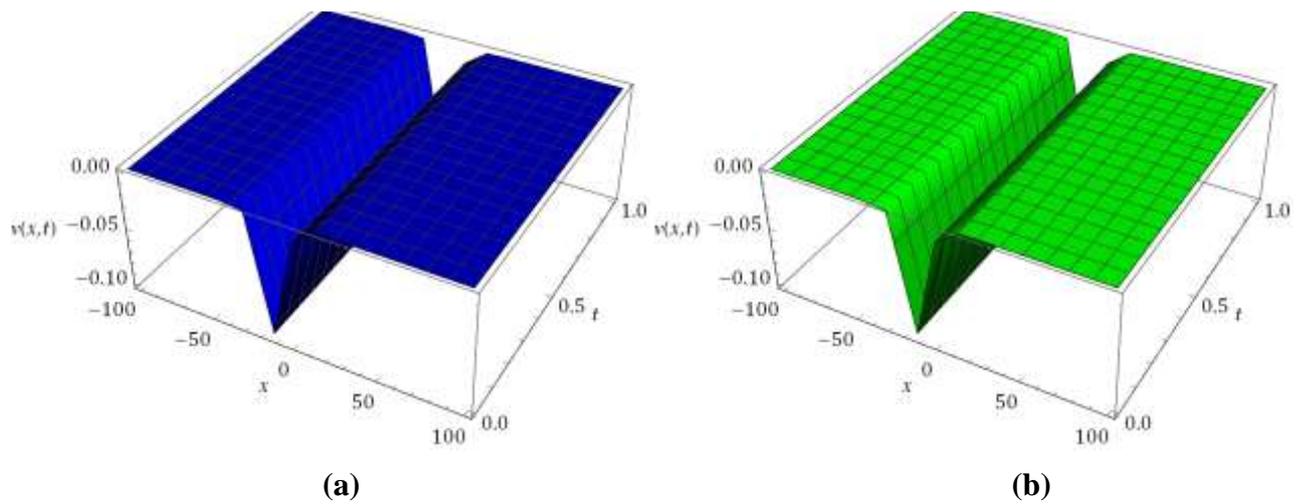

**(a)** **(b)**

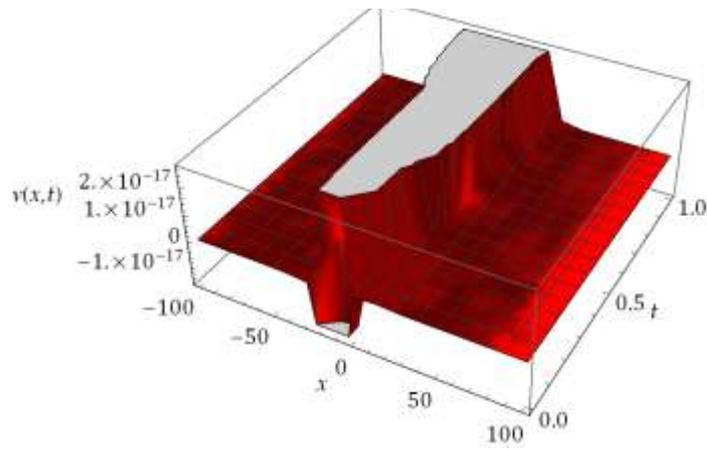

(c)

**Fig. 2.** Behaviour of **(a)** Approximate solution **(b)** Exact solution **(c)** Absolute error = $|v_{exa.} - v_{app.}|$ for Ex. 6.1 at $\omega = 0.005, \ell = 0.1, c = 10, n = 1, \alpha = 1$ and $\hbar = -1$.

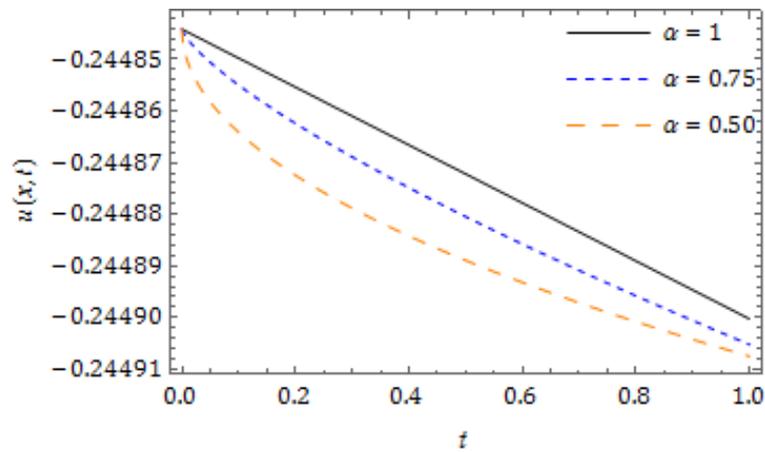

**Fig. 3**. Nature of $u(x,t)$ with $t$ for Ex. 6.1 at $\omega = 0.005, \hbar = -1, \ell = 0.1, c = 10, x = 1$ and $n = 1$ with diverse $\alpha$.

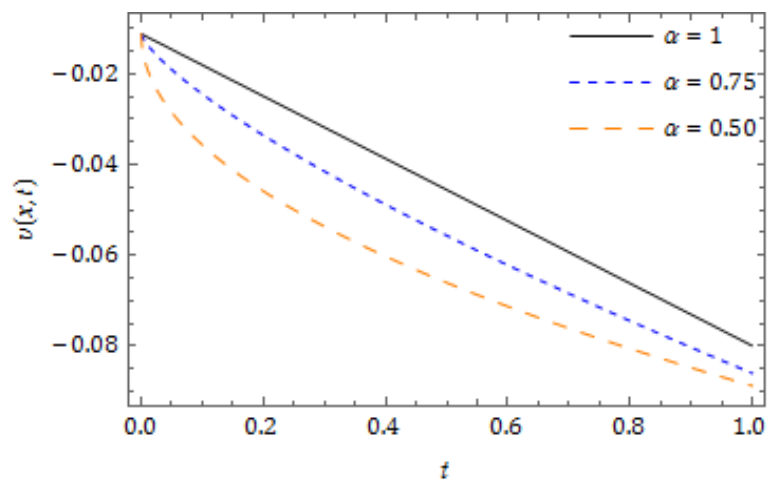

**Fig. 4**. Nature $v(x,t)$ with $t$ for Ex. 6.1 at $\omega = 0.005, \ell = 0.1, c = 10, \hbar = -1$  $x = 1$ and $n = 1$ with diverse $\alpha$.

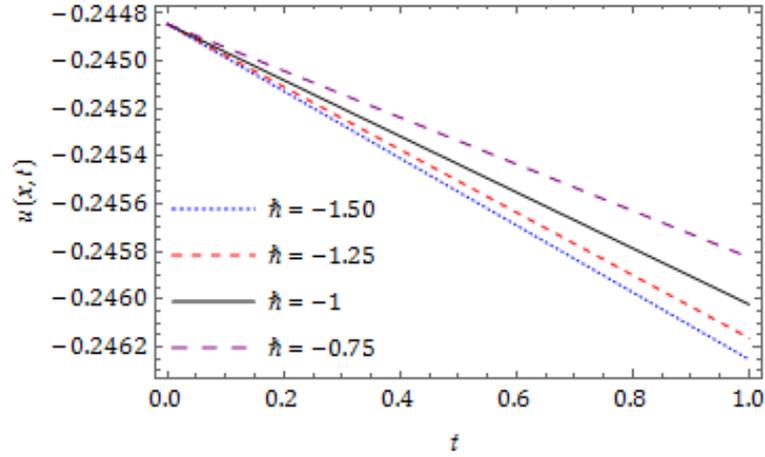

**Fig. 5**. Plot of $u(x,t)$ with diverse $\hbar$ when $\omega = 0.005, \ell = 0.1, c = 10, n = 5, \alpha = 1$ and $x = 1$ for Ex. 6.1.

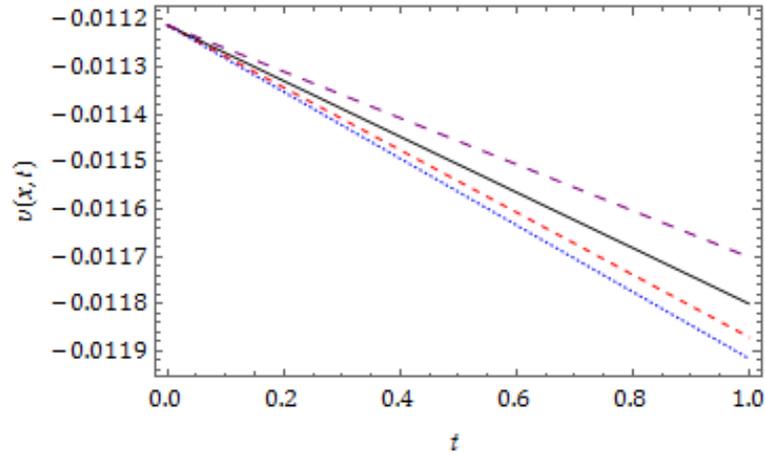

**Fig. 6**. Nature of $v(x,t)$ with distinct $\hbar$ at $\omega = 0.005, \ell = 0.1, c = 10, n = 1, \alpha = 1$ and $x = 1$ for Ex. 6.1 .

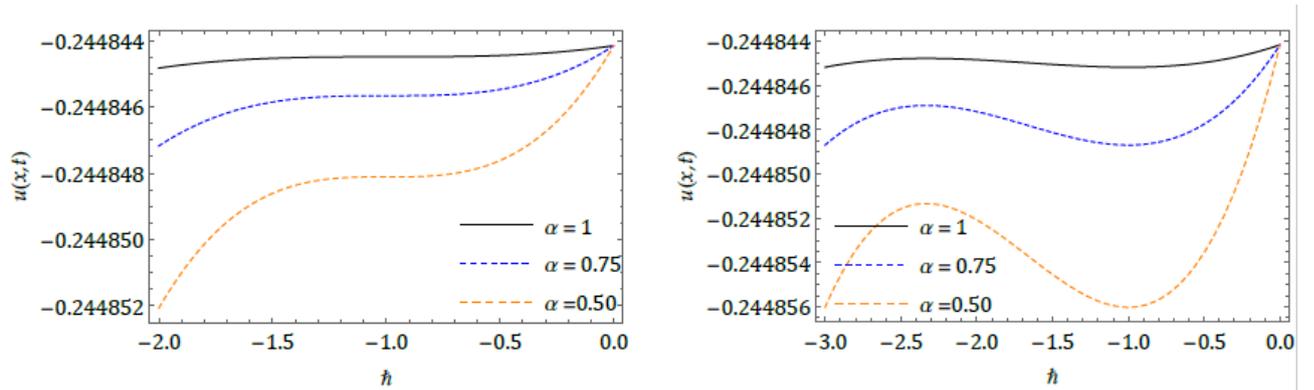

**Fig. 7**. $\hbar$-curves drown for $u(x,t)$ at $n = 1$ (left) and $n = 2$ (right) with different $\alpha$ when $\omega = 0.005, \ell = 0.1, c = 10, x = 1$ and $t = 0.01$ for Ex. 6.1 .

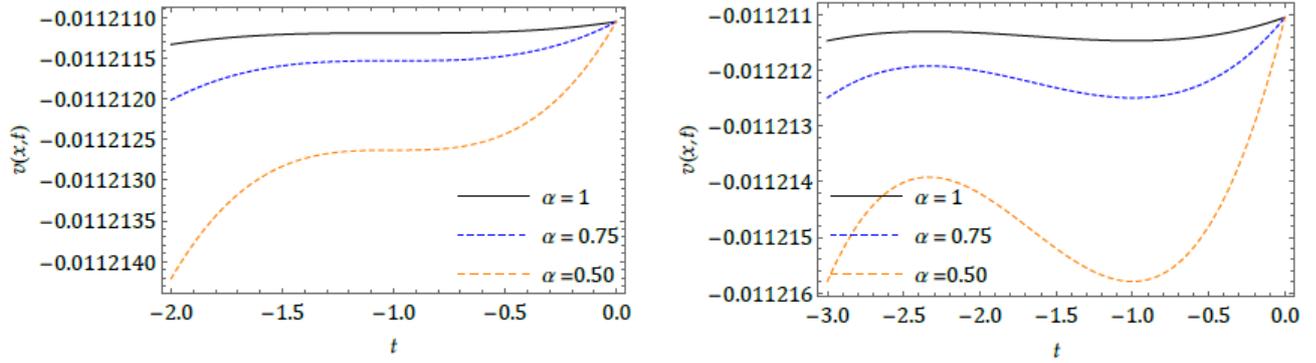

**Fig. 8.** ℏ-curve drown for $v(x,t)$ at $n = 1$ (left) and $n = 2$ (right) with different $\alpha$ when $\omega = 0.005, \ell = 0.1, c = 10, x = 1$ and $t = 0.01$ for Ex. 6.1.

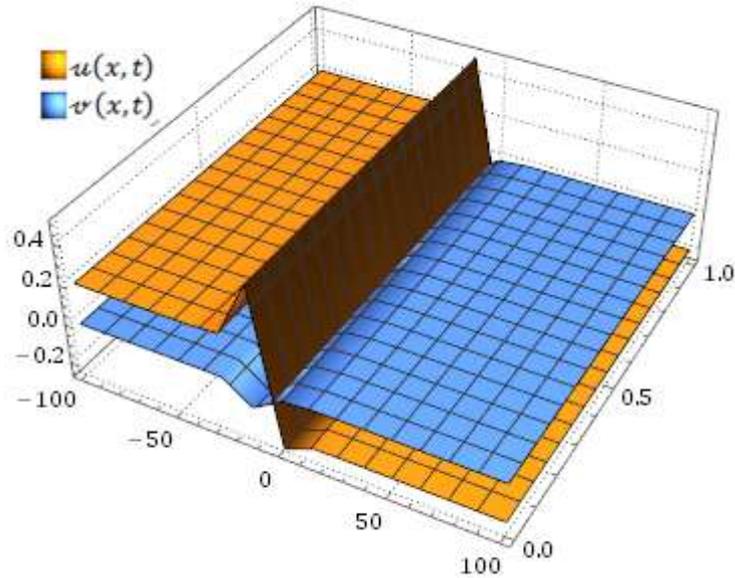

**Fig. 9.** Coupled surface of $u(x,t)$ and $v(x,t)$ for Ex. 6.1 at $\omega = 0.005, \ell = 0.1, c = 10, n = 1, \alpha = 1$ and $\hbar = -1$

**Table 1** Comparative study in terms of absolute error between ADM [17], VIM [18], CFRDTM [15] and $q$-HATM for the approximate solution $u(x,t)$ at $\omega = 0.005, \ell = 0.1, c = 10, n = 1, \hbar = -1$ and $\alpha = 1$ for Ex. 6.1.

| $(t,x)$ | $|u_{Exact} - u_{ADM}|$ | $|u_{Exact} - u_{VIM}|$ | $|u_{Exact} - u_{CRFDTM}|$ | $|u_{Exact} - u^{(3)}_{q-HATM}|$ |
|---|---|---|---|---|
| $(0.1, 0.1)$ | $8.16297 \times 10^{-7}$ | $6.35269 \times 10^{-5}$ | $5.55112 \times 10^{-17}$ | $5.55112 \times 10^{-17}$ |
| $(0.1, 0.3)$ | $7.64245 \times 10^{-7}$ | $1.90854 \times 10^{-4}$ | $5.55112 \times 10^{-17}$ | $5.55112 \times 10^{-17}$ |
| $(0.1, 0.5)$ | $7.16083 \times 10^{-7}$ | $3.18549 \times 10^{-4}$ | $5.55112 \times 10^{-16}$ | $5.55112 \times 10^{-16}$ |
| $(0.2, 0.1)$ | $3.26243 \times 10^{-6}$ | $6.18930 \times 10^{-5}$ | $5.55112 \times 10^{-16}$ | $5.55112 \times 10^{-16}$ |
| $(0.2, 0.3)$ | $3.05458 \times 10^{-6}$ | $1.85945 \times 10^{-4}$ | $1.11022 \times 10^{-16}$ | $1.11022 \times 10^{-16}$ |
| $(0.2, 0.5)$ | $2.86226 \times 10^{-6}$ | $3.10352 \times 10^{-4}$ | $7.77156 \times 10^{-16}$ | $7.77156 \times 10^{-16}$ |
| $(0.3, 0.1)$ | $7.33445 \times 10^{-6}$ | $6.03095 \times 10^{-5}$ | $0$ | $0$ |
| $(0.3, 0.3)$ | $6.86758 \times 10^{-6}$ | $1.81187 \times 10^{-4}$ | $1.66533 \times 10^{-16}$ | $1.66533 \times 10^{-16}$ |

| (0.3, 0.5) | $6.43557 \times 10^{-6}$ | $3.02408 \times 10^{-4}$ | $6.666134 \times 10^{-16}$ | $6.666134 \times 10^{-16}$ |
|---|---|---|---|---|
| (0.4, 0.1) | $1.30286 \times 10^{-5}$ | $5.87746 \times 10^{-5}$ | $5.55112 \times 10^{-17}$ | $5.55112 \times 10^{-17}$ |
| (0.4, 0.3) | $1.22000 \times 10^{-5}$ | $1.76574 \times 10^{-4}$ | $5.55112 \times 10^{-17}$ | $5.55112 \times 10^{-17}$ |
| (0.4, 0.5) | $1.14333 \times 10^{-5}$ | $2.94707 \times 10^{-4}$ | $5.55112 \times 10^{-16}$ | $5.55112 \times 10^{-16}$ |
| (0.5, 0.1) | $2.03415 \times 10^{-5}$ | $5.72867 \times 10^{-5}$ | 0 | 0 |
| (0.5, 0.3) | $1.90489 \times 10^{-5}$ | $1.72102 \times 10^{-4}$ | $1.11022 \times 10^{-16}$ | $1.11022 \times 10^{-16}$ |
| (0.5, 0.5) | $1.78528 \times 10^{-5}$ | $2.87241 \times 10^{-4}$ | $6.10623 \times 10^{-16}$ | $6.10623 \times 10^{-16}$ |

**Table 2** Comparative study in terms of absolute error between ADM [17], VIM [18], CFRDTM [15] and $q$-HATM for the approximate solution $v(x,t)$ at $\omega = 0.005, \ell = 0.1, c = 10, \hbar = -1, n = 1, \hbar = -1$ and $\alpha = 1$ for Ex. 6.1.

| $(t,x)$ | $\|v_{Exact} - v_{ADM}\|$ | $\|v_{Exact} - v_{VIM}\|$ | $\|v_{Exact} - v_{CRFDTM}\|$ | $\|v_{Exact} - v^{(3)}_{q-HATM}\|$ |
|---|---|---|---|---|
| (0.1, 0.1) | $5.88676 \times 10^{-5}$ | $1.65942 \times 10^{-5}$ | $3.46945 \times 10^{-18}$ | $3.46945 \times 10^{-18}$ |
| (0.1, 0.3) | $5.56914 \times 10^{-5}$ | $4.98691 \times 10^{-5}$ | $5.55112 \times 10^{-17}$ | $5.55112 \times 10^{-17}$ |
| (0.1, 0.5) | $5.27169 \times 10^{-5}$ | $8.32598 \times 10^{-5}$ | $5.55112 \times 10^{-16}$ | $5.55112 \times 10^{-16}$ |
| (0.2, 0.1) | $1.18213 \times 10^{-4}$ | $1.06813 \times 10^{-5}$ | $6.93889 \times 10^{-18}$ | $6.93889 \times 10^{-18}$ |
| (0.2, 0.3) | $1.11833 \times 10^{-4}$ | $4.83269 \times 10^{-5}$ | $5.55112 \times 10^{-17}$ | $5.55112 \times 10^{-17}$ |
| (0.2, 0.5) | $1.05858 \times 10^{-4}$ | $8.06837 \times 10^{-5}$ | $5.55112 \times 10^{-16}$ | $5.55112 \times 10^{-16}$ |
| (0.3, 0.1) | $1.78041 \times 10^{-4}$ | $1.55880 \times 10^{-5}$ | $6.93889 \times 10^{-18}$ | $6.93889 \times 10^{-18}$ |
| (0.3, 0.3) | $1.68429 \times 10^{-4}$ | $4.68440 \times 10^{-5}$ | $5.55112 \times 10^{-17}$ | $5.55112 \times 10^{-17}$ |
| (0.3, 0.5) | $1.59428 \times 10^{-4}$ | $7.82068 \times 10^{-5}$ | $5.55112 \times 10^{-16}$ | $5.55112 \times 10^{-16}$ |
| (0.4, 0.1) | $2.38356 \times 10^{-4}$ | $1.51135 \times 10^{-5}$ | $5.20417 \times 10^{-18}$ | $5.20417 \times 10^{-18}$ |
| (0.4, 0.3) | $2.25483 \times 10^{-4}$ | $4.54174 \times 10^{-5}$ | $5.55112 \times 10^{-17}$ | $5.55112 \times 10^{-17}$ |
| (0.4, 0.5) | $2.13430 \times 10^{-4}$ | $7.58243 \times 10^{-5}$ | $5.55112 \times 10^{-16}$ | $5.55112 \times 10^{-16}$ |
| (0.5, 0.1) | $2.99162 \times 10^{-4}$ | $1.46569 \times 10^{-5}$ | $1.73472 \times 10^{-18}$ | $1.73472 \times 10^{-18}$ |
| (0.5, 0.3) | $2.83001 \times 10^{-4}$ | $4.40448 \times 10^{-5}$ | $5.55112 \times 10^{-17}$ | $5.55112 \times 10^{-17}$ |
| (0.5, 0.5) | $2.67868 \times 10^{-4}$ | $7.35317 \times 10^{-5}$ | $5.55112 \times 10^{-16}$ | $5.55112 \times 10^{-16}$ |

**Example 6.2**.

Consider the approximate long wave (ALW) equations with arbitrary order [15, 17, 22]:

$$\begin{cases} D_t^\alpha u = -u\frac{\partial u}{\partial x} - \frac{\partial v}{\partial x} - \frac{1}{2}\frac{\partial^2 v}{\partial x^2}, \\ D_t^\alpha v = -u\frac{\partial v}{\partial x} - v\frac{\partial u}{\partial x} + \frac{1}{2}\frac{\partial^2 v}{\partial x^2}, \end{cases} \quad 0 < \alpha \leq 1, \quad (34)$$

with initial conditions

$$u(x, 0) = \omega - \ell \coth[\ell(x + c)], \quad v(x, 0) = -\ell^2 \operatorname{csch}^2[\ell(x + c)]. \quad (35)$$

The exact solution for classical oreder ALW equations is

$$u(x, t) = \omega - \ell \coth[\ell(x + c - \omega t)], \quad v(x, t) = -\ell^2 \operatorname{csch}^2[\ell(x + c - \omega t)].$$

By performing LT on Eq. (34) and using Eq. (35), we have

$$L[u(x,t)] - \frac{\omega - \ell coth[\ell(x+c)]}{s} + \frac{1}{s^\alpha} L\left\{ u\frac{\partial u}{\partial x} + \frac{\partial v}{\partial x} + \frac{1}{2}\frac{\partial^2 u}{\partial x^2} \right\} = 0,$$

$$L[v(x,t)] - \frac{-\ell^2 csch^2[\ell(x+c)]}{s} + \frac{1}{s^\alpha} L\left\{ u\frac{\partial v}{\partial x} + v\frac{\partial u}{\partial x} - \frac{1}{2}\frac{\partial^2 v}{\partial x^2} \right\} = 0. \tag{36}$$

Define the non-linear operators as

$$N^1[\varphi_1(x,t;q), \varphi_2(x,t;q)] = L[\varphi_1(x,t;q)] - \frac{\omega - \ell coth[\ell(x+c)]}{s}$$
$$+ \frac{1}{s^\alpha} L\left\{ \varphi_1(x,t;q)\frac{\partial \varphi_1(x,t;q)}{\partial x} + \frac{\partial \varphi_2(x,t;q)}{\partial y} + \frac{1}{2}\frac{\partial^2 \varphi_1(x,t;q)}{\partial x^2} \right\},$$

$$N^2[\varphi_1(x,t;q), \varphi_2(x,t;q)] = L[\varphi_2(x,t;q)] - \frac{-\ell^2 csch^2[\ell(x+c)]}{s} + \frac{1}{s^\alpha} L\{\varphi_1(x,t;q)\frac{\partial \varphi_2(x,t;q)}{\partial x}$$
$$+ \varphi_2(x,t;q)\frac{\partial \varphi_1(x,t;q)}{\partial x} - \frac{1}{2}\frac{\partial^2 \varphi_2(x,t;q)}{\partial x^2}\}. \tag{37}$$

Now, for $H(x,t) = 1$, the $m$-$th$ order deformation equation is presented as follows

$$L[u_m(x,t) - k_m u_{m-1}(x,t)] = \hbar \Re_{1,m}[\vec{u}_{m-1}, \vec{v}_{m-1}],$$
$$L[v_m(x,t) - k_m v_{m-1}(x,t)] = \hbar \Re_{2,m}[\vec{u}_{m-1}, \vec{v}_{m-1}]. \tag{38}$$

where

$$\Re_{1,m}[\vec{u}_{m-1}, \vec{v}_{m-1}] = L[u_{m-1}(x,t)] - \left(1 - \frac{k_m}{n}\right)\frac{\omega - \ell coth[\ell(x+c)]}{s}$$
$$+ \frac{1}{s^\alpha} L\left\{ \sum_{i=0}^{m-1} u_i \frac{\partial u_{m-1-i}}{\partial x} + \frac{\partial v_{m-1}}{\partial x} + \frac{1}{2}\frac{\partial^2 u_{m-1}}{\partial x^2} \right\},$$

$$\Re_{2,m}[\vec{u}_{m-1}, \vec{v}_{m-1}] = L[v_{m-1}(x,t)] - \left(1 - \frac{k_m}{n}\right)\frac{-\ell^2 csch^2[\ell(x+c)]}{s}$$
$$+ \frac{1}{s^\alpha} L\left\{ \sum_{i=0}^{m-1} u_i \frac{\partial v_{m-1-i}}{\partial x} + \sum_{i=0}^{m-1} v_i \frac{\partial u_{m-1-i}}{\partial x} - \frac{1}{2}\frac{\partial^2 v_{m-1}}{\partial x^2} \right\}. \tag{39}$$

By employing inverse LT on Eq. (38), we obtained

$$u_m(x,t) = k_m u_{m-1}(x,t) + \hbar L^{-1}\{\Re_{1,m}[\vec{u}_{m-1}, \vec{v}_{m-1}]\},$$
$$v_m(x,t) = k_m v_{m-1}(x,t) + \hbar L^{-1}\{\Re_{2,m}[\vec{u}_{m-1}, \vec{v}_{m-1}]\}. \tag{40}$$

On simplification, we obtain

$$u_0(x,t) = \omega - \ell coth[\ell(x+c)], \quad v_0(x,t) = -\ell^2 csch^2[\ell(x+c)],$$

$$u_1(x,t) = \frac{\hbar \ell^2 \omega\, csch^2[\ell(x+c)] t^\alpha}{\Gamma[\alpha+1]}, \quad v_1(x,t) = \frac{2\hbar \ell^3 \omega\, coth[\ell(x+c)]\, csch^2[\ell(x+c)] t^\alpha}{\Gamma[\alpha+1]},$$

$$u_2(x,t) = \frac{(n+\hbar)\hbar \ell^2 \omega\, csch^2[\ell(x+c)] t^\alpha}{\Gamma[\alpha+1]} - \frac{2\hbar^2 \ell^3 \omega^2\, coth[\ell(x+c)] csch^2[\ell(x+c)] t^{2\alpha}}{\Gamma[2\alpha+1]},$$

$$v_2(x,t) = \frac{2(n+\hbar)\hbar \ell^3 \omega\, coth[\ell(x+c)]\, csch^2[\ell(x+c)] t^\alpha}{\Gamma[\alpha+1]} - \frac{2\hbar^2 \ell^4 \omega^2 (2+\cosh[2\ell(x+c)])\, csch^4[\ell(x+c)] t^{2\alpha}}{\Gamma[2\alpha+1]},$$

$$u_3(x,t) = \frac{(\hbar+n)^2 \hbar \ell^2 \omega\, csch^2[\ell(x+c)]\, t^\alpha}{\Gamma[1+\alpha]} - \frac{4(\hbar+n)\hbar^2 \ell^3 \omega^2\, coth[\ell(x+c)]\, csch^2[\ell(x+c)]\, t^{2\alpha}}{\Gamma[1+2\alpha]}$$
$$+ \frac{\hbar^3 \ell^4 t^{3\alpha} \omega^2\, csch^5[\ell(x+c)]}{\Gamma[\alpha+1]^2\, \Gamma[3\alpha+1]}(-2k\, cosh[\ell(x+c)]\, \Gamma[1+2\alpha]$$

$$+\Gamma[1+\alpha]^2\big(4\ell\,\cosh[\ell(x+c)]+\omega(3\,\sinh[\ell(x+c)]+\sinh[3\ell(x+c)])\big),$$

$$v_3(x,t) = \frac{2(\hbar+n)^2\hbar\ell^3\omega\,\coth[\ell(x+c)]\,csch^2[\ell(x+c)]\,t^\alpha}{\Gamma[\alpha+1]} - \frac{4(\hbar+n)\hbar^2\ell^5 t^{2\beta}\omega(2+\cosh[2\ell(x+c)])\,csch^4[\ell(x+c)]}{\Gamma[2\alpha+1]}$$

$$+\frac{\hbar^3\ell^5 t^{3\alpha}\omega^2\,csch^6[\ell(x+c)]}{\Gamma[\alpha+1]^2\,\Gamma[3\alpha+1]}(-2l(3+2\cosh[2\ell(x+c)])\Gamma[2\alpha+1]$$

$$+\Gamma[\alpha+1]^2(12\ell+8\ell\cosh[2\ell(x+c)]+10\omega\sinh[2\ell(x+c)]+\omega\sinh[4\ell(x+c)]),$$

$\vdots$

Following in the same procedure, the remaining terms can be obtained. Eventually, for Eq. (34) the $q$-HATM solutions is presented as follows

$$u(x,t) = u_0(x,t) + \sum_{m=1}^{\infty} u_m(x,t)\left(\frac{1}{n}\right)^m,$$
$$v(x,t) = v_0(x,t) + \sum_{m=1}^{\infty} v_m(x,t)\left(\frac{1}{n}\right)^m. \tag{41}$$

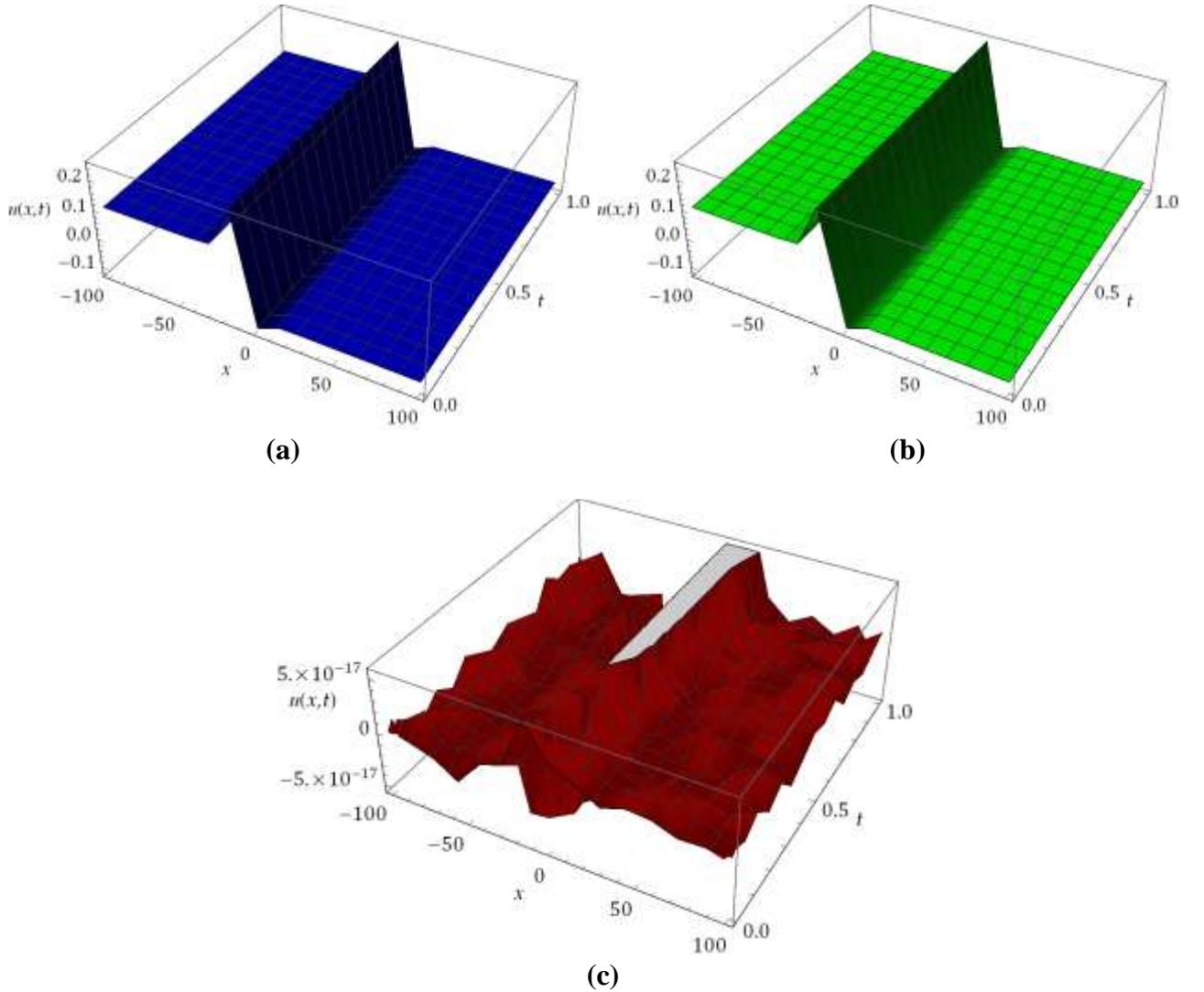

**Fig. 10**. Behaviour of **(a)** Approximate solution **(b)** Exact solution **(c)** Absolute error = $|u_{exa.} - u_{app.}|$ at $\omega = 0.005, \ell = 0.1, c = 10, n = 1, \alpha = 1$ and $\hbar = -1$ for Ex. 6.2.

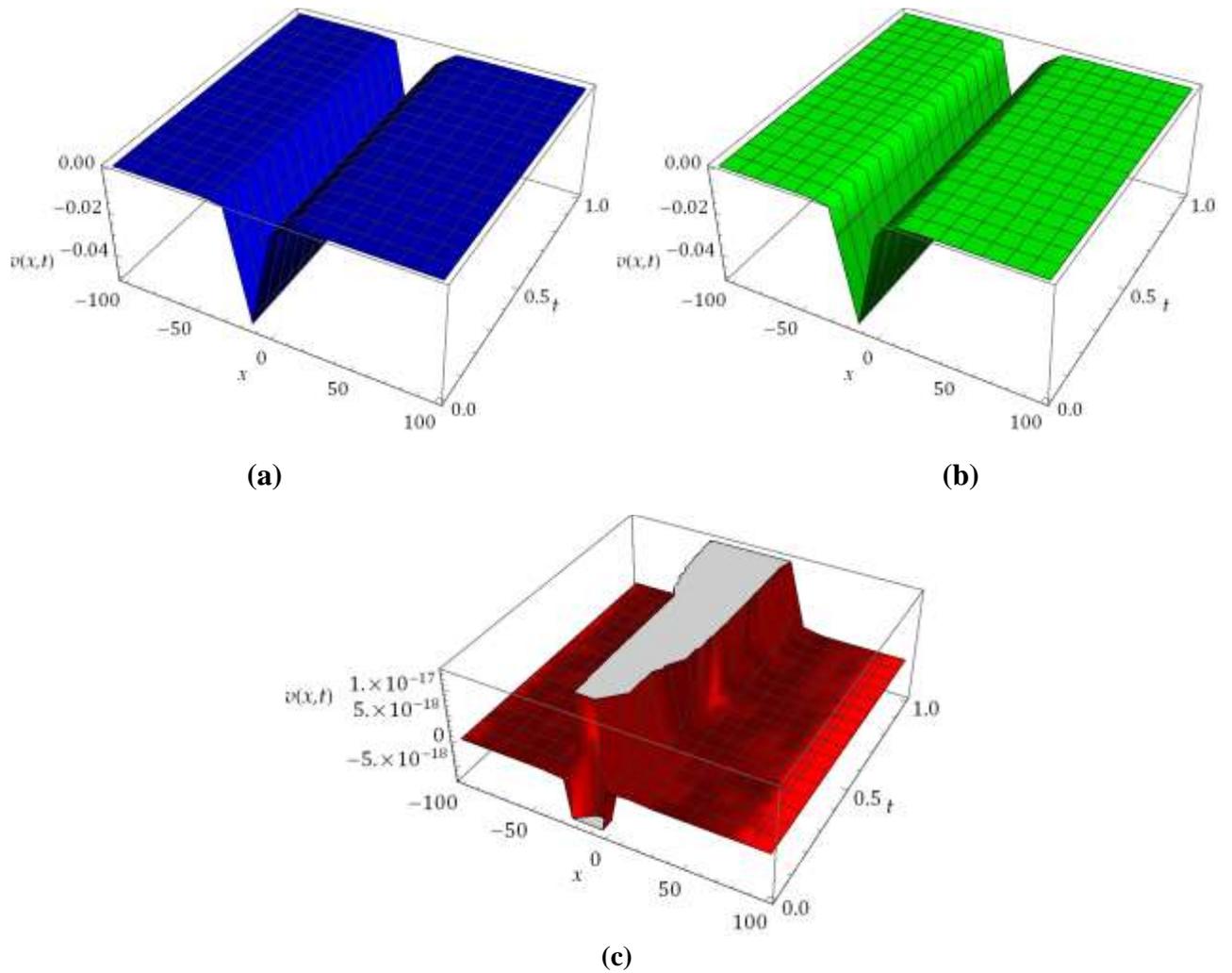

**Fig. 11**. Behaviour of **(a)** Approximate solution **(b)** Exact solution **(c)** Absolute error = $|v_{exa.} - v_{app.}|$ at $\omega = 0.005, \ell = 0.1, c = 10, n = 1, \alpha = 1$ and $\hbar = -1$ for Ex. 6.2.

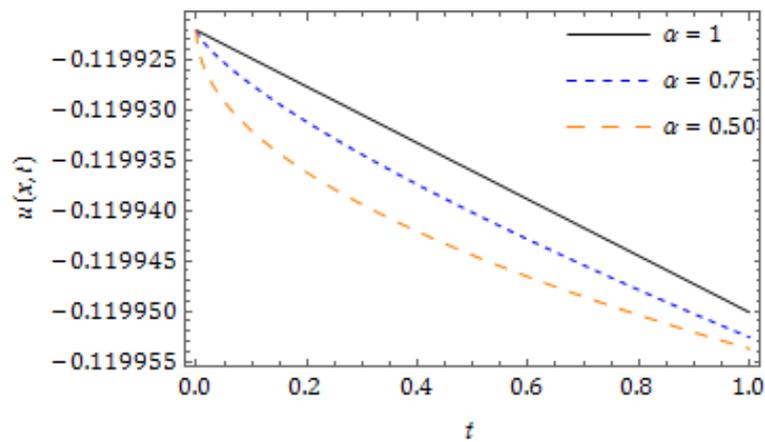

**Fig. 12**. Plot of $u(x,t)$ with respect to $t$ at diverse values of $\alpha$ when $\omega = 0.005, \ell = 0.1, c = 10, n = 1, \hbar = -1$ and $x = 1$ for Ex. 6.2.

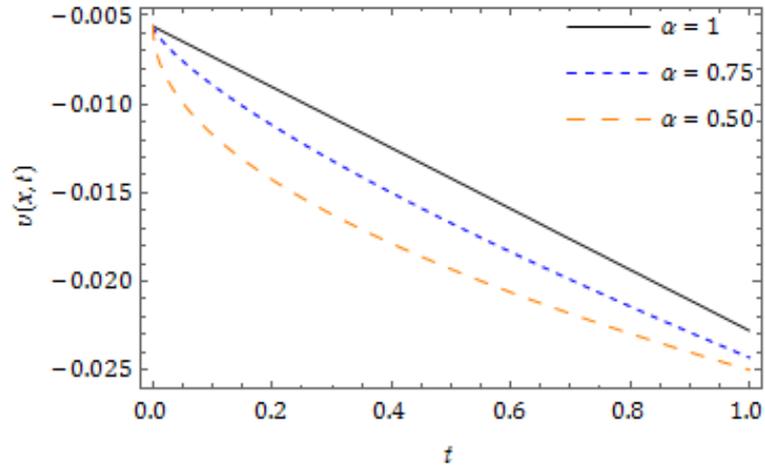

**Fig. 13**. Response of $v(x,t)$ with $t$ at diverse $\alpha$ when $\omega = 0.005, \ell = 0.1, c = 10, n = 1, \hbar = -1$ and $x = 1$ for Ex. 6.2.

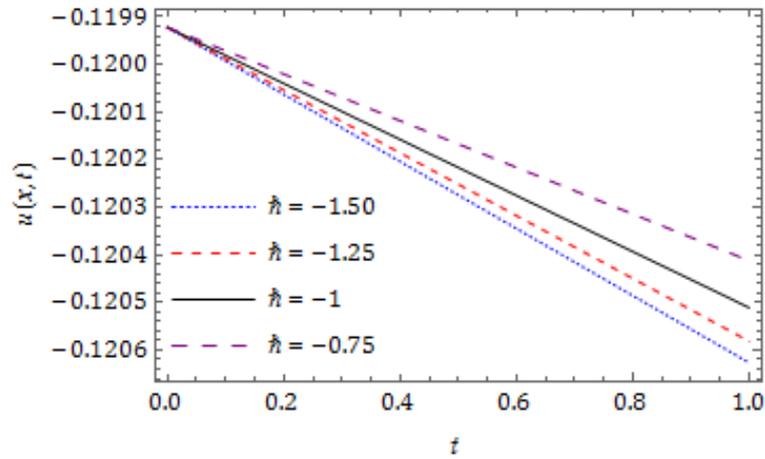

**Fig. 14**. Nature of $u(x,t)$ with distinct $\hbar$ for Ex. 6.2 at $\omega = 0.005, \ell = 0.1, c = 10, n = 5, \alpha = 1$ and $x = 1$.

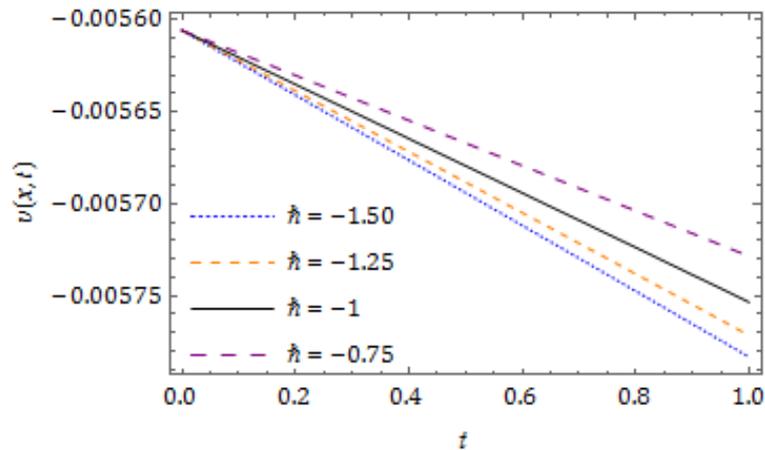

**Fig. 15**. Nature of $v(x,t)$ with distinct $\hbar$ for Ex. 6.2 at $\ell = 0.1, c = 10, n = 3, \alpha = 1$ and $x = 1$.

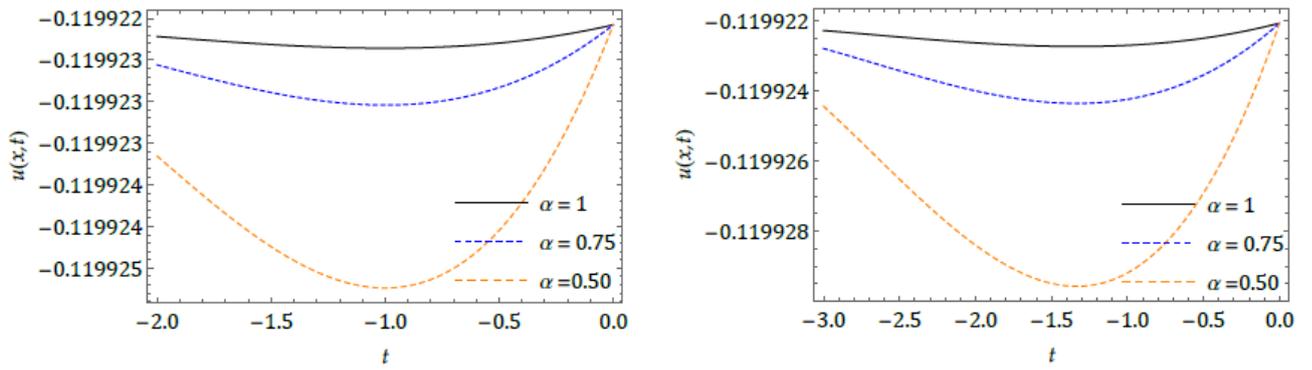

**Fig. 16.** $\hbar$-curves drown for $u(x,t)$ at $n=1$ (left) and $n=2$ (right) for Ex. 6.2 at $\omega=0.005, \ell=0.1, c=10, x=1$ and $t=0.01$ with diverse $\alpha$.

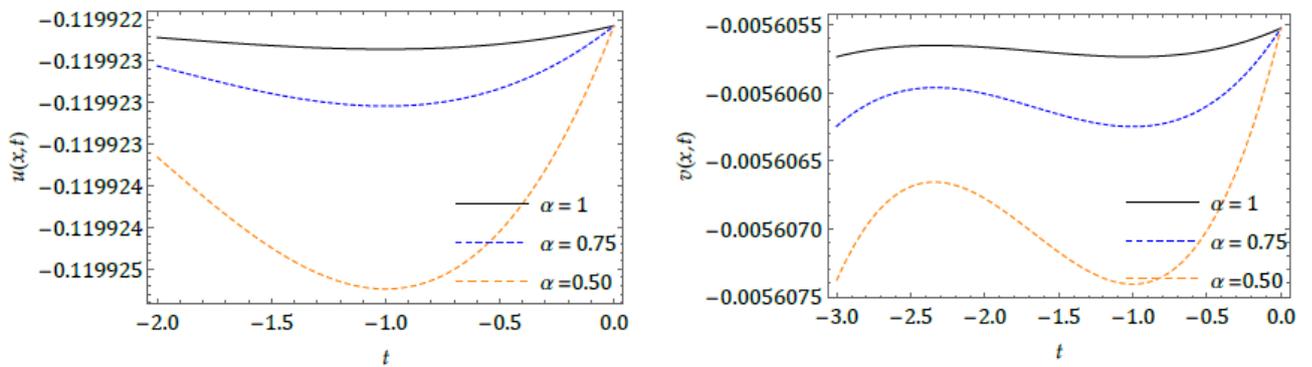

**Fig. 17.** $\hbar$-curves drown for $v(x,t)$ at $n=1$ (left) and $n=2$ (right) with different $\alpha$ for Ex. 6.1 when $\omega=0.005, \ell=0.1, c=10, x=1$ and $t=0.01$

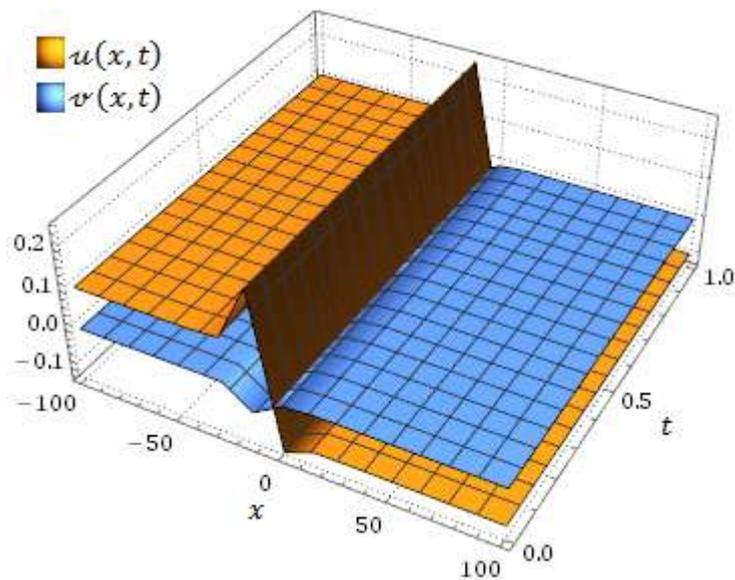

**Fig. 18.** Nature of coupled $q$-HATM solutions $u(x,t)$ and $v(x,t)$ for Ex. 6.2 when at $n=1, \alpha=1, \hbar=-1, \omega=0.005, \ell=0.1$ and $c=10$.

**Table 3** Comparative study between ADM [17], VIM [18], LADM [19], CRFDTM [15] and $q$-HATM for the approximate solution $u(x,t)$ at $\omega = 0.005, \ell = 0.1\ c = 10, \alpha = 1$ and $\hbar = -1$ for Ex. 6.2.

| $(x,t)$ | $\|u_{Exact} - u_{ADM}\|$ | $\|u_{Exact} - u_{VIM}\|$ | $\|u_{Exact} - u_{LADM}\|$ | $\|u_{Exact} - u_{CRFDTM}\|$ | $\|u_{Exact} - u^{(3)}_{q-HATM}\|$ |
|---|---|---|---|---|---|
| (0.1, 0.1) | $8.02989 \times 10^{-6}$ | $1.23033 \times 10^{-4}$ | $7.10000 \times 10^{-9}$ | $2.77556 \times 10^{-17}$ | $2.77556 \times 10^{-17}$ |
| (0.1, 0.3) | $7.38281 \times 10^{-6}$ | $3.69597 \times 10^{-4}$ | $6.50000 \times 10^{-9}$ | $2.77556 \times 10^{-17}$ | $2.77556 \times 10^{-17}$ |
| (0.1, 0.5) | $6.79923 \times 10^{-6}$ | $4.92780 \times 10^{-4}$ | $5.90000 \times 10^{-9}$ | $3.33067 \times 10^{-16}$ | $3.33067 \times 10^{-16}$ |
| (0.2, 0.1) | $3.23228 \times 10^{-5}$ | $1.69274 \times 10^{-5}$ | $2.82000 \times 10^{-8}$ | $2.77556 \times 10^{-17}$ | $2.77556 \times 10^{-17}$ |
| (0.2, 0.3) | $297172 \times 10^{-5}$ | $1.89210 \times 10^{-4}$ | $2.59000 \times 10^{-8}$ | $4.16334 \times 10^{-17}$ | $4.16334 \times 10^{-17}$ |
| (0.2, 0.5) | $2.73673 \times 10^{-5}$ | $1.55176 \times 10^{-4}$ | $2.41000 \times 10^{-8}$ | $3.60822 \times 10^{-17}$ | $3.60822 \times 10^{-17}$ |
| (0.3, 0.1) | $7.32051 \times 10^{-5}$ | $1.12345 \times 10^{-5}$ | $6.33670 \times 10^{-8}$ | $1.38778 \times 10^{-17}$ | $1.38778 \times 10^{-17}$ |
| (0.3, 0.3) | $6.73006 \times 10^{-5}$ | $6.55176 \times 10^{-5}$ | $5.85000 \times 10^{-8}$ | $2.77556 \times 10^{-17}$ | $2.77556 \times 10^{-17}$ |
| (0.3, 0.5) | $6.19760 \times 10^{-5}$ | $2.12346 \times 10^{-5}$ | $5.40000 \times 10^{-8}$ | $3.19189 \times 10^{-16}$ | $3.19189 \times 10^{-16}$ |
| (0.4, 0.1) | $1.31032 \times 10^{-4}$ | $7.36513 \times 10^{-5}$ | $1.12400 \times 10^{-7}$ | $1.38778 \times 10^{-17}$ | $1.38778 \times 10^{-17}$ |
| (0.4, 0.3) | $1.20455 \times 10^{-4}$ | $9.5016 \times 10^{-5}$ | $1.03900 \times 10^{-7}$ | $2.77556 \times 10^{-17}$ | $2.77556 \times 10^{-17}$ |
| (0.4, 0.5) | $1.10919 \times 10^{-4}$ | $8.23160 \times 10^{-4}$ | $9.61000 \times 10^{-8}$ | $3.19189 \times 10^{-16}$ | $3.19189 \times 10^{-16}$ |
| (0.5, 0.1) | $2.06186 \times 10^{-4}$ | $5.55176 \times 10^{-5}$ | $1.75500 \times 10^{-7}$ | 0 | 0 |
| (0.5, 0.3) | $1.89528 \times 10^{-4}$ | $3.21715 \times 10^{-6}$ | $1.62200 \times 10^{-7}$ | $5.55112 \times 10^{-17}$ | $5.55112 \times 10^{-17}$ |
| (0.5, 0.5) | $1.74510 \times 10^{-4}$ | $2.00176 \times 10^{-5}$ | $1.5010 \times 10^{-7}$ | $3.19189 \times 10^{-16}$ | $3.19189 \times 10^{-16}$ |

**Table 4** Comparative study between ADM [17], VIM [18], LADM [19], CRFDTM [15] and $q$-HATM for the approximate solution $v(x,t)$ at $\omega = 0.005, \ell = 0.1\ c = 10, \alpha = 1$ and $\hbar = -1$ for Ex. 6.2.

| $(x,t)$ | $\|v_{Exact} - v_{ADM}\|$ | $\|v_{Exact} - v_{VIM}\|$ | $\|v_{Exact} - v_{LADM}\|$ | $\|v_{Exact} - v_{CRFDTM}\|$ | $\|v_{Exact} - v^{(3)}_{q-HATM}\|$ |
|---|---|---|---|---|---|
| (0.1, 0.1) | $4.81902 \times 10^{-4}$ | $1.23033 \times 10^{-4}$ | $9.5512 \times 10^{-10}$ | $1.73472 \times 10^{-18}$ | $1.73472 \times 10^{-18}$ |
| (0.1, 0.3) | $4.50818 \times 10^{-4}$ | $1.7600 \times 10^{-4}$ | $8.0600 \times 10^{-10}$ | $2.60209 \times 10^{-17}$ | $2.60209 \times 10^{-17}$ |
| (0.1, 0.5) | $4.22221 \times 10^{-4}$ | $2.69597 \times 10^{-4}$ | $6.7700 \times 10^{-10}$ | $1.80411 \times 10^{-16}$ | $1.80411 \times 10^{-16}$ |
| (0.2, 0.1) | $9.76644 \times 10^{-4}$ | $2.69597 \times 10^{-4}$ | $3.8210 \times 10^{-9}$ | $3.46945 \times 10^{-18}$ | $3.46945 \times 10^{-18}$ |
| (0.2, 0.3) | $9.13502 \times 10^{-4}$ | $2.69597 \times 10^{-4}$ | $3.224 \times 10^{-9}$ | $2.34188 \times 10^{-17}$ | $2.34188 \times 10^{-17}$ |
| (0.2, 0.5) | $8.55426 \times 10^{-4}$ | $2.69597 \times 10^{-4}$ | $2.7060 \times 10^{-9}$ | $1.73472 \times 10^{-16}$ | $1.73472 \times 10^{-16}$ |
| (0.3, 0.1) | $1.48482 \times 10^{-3}$ | $2.69597 \times 10^{-4}$ | $8.597 \times 10^{-9}$ | $3.46945 \times 10^{-18}$ | $3.46945 \times 10^{-18}$ |
| (0.3, 0.3) | $1.38858 \times 10^{-3}$ | $2.69597 \times 10^{-4}$ | $7.252 \times 10^{-9}$ | $1.99493 \times 10^{-17}$ | $1.99493 \times 10^{-17}$ |
| (0.3, 0.5) | $1.30009 \times 10^{-3}$ | $2.69597 \times 10^{-4}$ | $6.0910 \times 10^{-9}$ | $1.61329 \times 10^{-16}$ | $1.61329 \times 10^{-16}$ |
| (0.4, 0.1) | $2.00705 \times 10^{-3}$ | $2.69597 \times 10^{-4}$ | $1.5284 \times 10^{-8}$ | $2.60209 \times 10^{-18}$ | $2.60209 \times 10^{-18}$ |
| (0.4, 0.3) | $1.87661 \times 10^{-3}$ | $2.69597 \times 10^{-4}$ | $1.2893 \times 10^{-8}$ | $1.73472 \times 10^{-17}$ | $1.73472 \times 10^{-17}$ |
| (0.4, 0.5) | $1.75670 \times 10^{-3}$ | $2.69597 \times 10^{-4}$ | $1.0827 \times 10^{-8}$ | $1.52656 \times 10^{-16}$ | $1.52656 \times 10^{-16}$ |
| (0.5, 0.1) | $2.54396 \times 10^{-3}$ | $2.69597 \times 10^{-4}$ | $2.3880 \times 10^{-8}$ | $8.67362 \times 10^{-19}$ | $8.67362 \times 10^{-19}$ |
| (0.5, 0.3) | $2.37815 \times 10^{-3}$ | $2.69597 \times 10^{-4}$ | $2.0144 \times 10^{-8}$ | $2.08167 \times 10^{-17}$ | $2.08167 \times 10^{-17}$ |
| (0.5, 0.5) | $2.22578 \times 10^{-3}$ | $2.69597 \times 10^{-4}$ | $1.6916 \times 10^{-8}$ | $1.43982 \times 10^{-16}$ | $1.43982 \times 10^{-16}$ |

## 7. Numerical results and discussion

Here, the numerical simulation has been conducted in order to prove whether the future algorithm lead to greater accuracy. We can see that from obtained results, the future scheme gives remarkable exactness in comparison to the technique presented in the literature [15,18, 19], which is cited for both the examples in the Tables 1- 4. Figs. 1 and 2 explore the comparison of obtained solutions with exact solutions and absolute error for Ex. 6.1. Figs. 3 and 4 are the response of obtained solutions for FMB equation with diverse Brownian motion and standard motion ($\alpha = 1$). Figs. 5 and 6 depict the $q$-HATM solutions for distinct $\hslash$, which aids us to adjust and control the convergence region. The Figs. 7-8 presents the importance of asymptotic parameter $n$ with $\hbar$ in $q$-HATM solution.

Moreover, Figs. 10 and 11 show the nature of obtained solutions in comparison with exact solutions for Ex. 6.2, in particular Figs. 10 (c) and 11 (c) revel the exactness of the obtained solution in absolute error. Figs. 12 and 13 explore the behaviour of obtained solution for distinct $\alpha$ i.e., $\alpha = 1, 0.75$ and $0.50$. Figs. 14 and 15 cites the nature of obtained solutions for distinct $\hbar$, and these aids us to control the convergence region. Finally, Figs. 16-17 signify the $\hbar$-curves and the horizontal line illustrate the region of convergence for FALW equation. The coupled surface of the MB and ALW equations consider in Ex. 6.1 and Ex.6.2 are respectively shown in Figs. 9 and 18, which help understand the nature of coupled equations.

## 8. Conclusion:

In the present work, the $q$-HATM is employed advantageously to find the solution for coupled modified Boussinesq and approximate long wave equations of fractional order. Two examples are considered in order to illustrate and validate the efficiency of the considered algorithm. The convergence and error analysis have been offered in the present frame work to show the consistency and applicability. The numerical simulation has been conducted for both considered fractional coupled systems in terms of absolute error. From the cited tables and plots, we can see that the proposed technique is effective and more exact as related to other methods and contains the results of CFRDTM as a special case ($n = 1$ and $\hbar$). Moreover, the algorithm controls and manipulates the series solution and which quickly converges to analytical solution in a small admissible domain. Hence, we can concluded that the proposed algorithm is very powerful and well organized to study the coupled system arise in physical phenomena both fractional and integer order derivative by analytically and numerically to describe the real world problems in a systematic and better manner.